\documentclass[11pt]{article}

\usepackage{latexsym}
\usepackage{amsmath}
\usepackage{amsfonts}
\usepackage{graphicx,epsfig}
\usepackage{color}

\newcommand{\bitem}{\begin{itemize}}
\newcommand{\eitem}{\end{itemize}}
\newcommand{\beq}{\begin{equation}}
\newcommand{\eeq}{\end{equation}}
\newcommand{\beqn}{\begin{eqnarray*}}
\newcommand{\eeqn}{\end{eqnarray*}}
\newcommand{\goto}{\rightarrow}
\newcommand{\cC}{{\cal C}}

\newcommand{\cP}{{\cal P}}

\newcommand{\cS}{{\cal S}}

\newcommand{\cK}{{\cal K}}
\newcommand{\bR}{{\bf R}}
\newcommand{\cE}{{\cal E}}

\newcommand{\Sp}{\mbox{supp} \,}

\def\t{\tilde}
\newcommand{\ip}[2]{\left\langle#1,#2\right\rangle}

\newcommand{\absip}[2]{| \langle#1,#2\rangle |}
\newcommand{\norm}[1]{\|#1\|}
\newcommand{\qed}{$\Box$}

\def\diag{{\text{\rm diag}}}

\def\cH{\mathcal{H}}

\def\N{\mathbb{N}}

\def\ZZ{\mathbb{Z}}
\def\bZ{\mathbb{Z}}

\def\R{\mathbb{R}}
\def\RR{\mathbb{R}}
\def\bR{\mathbb{R}}

\def\cZ{\mathcal{Z}}

\def\t{\tilde}

\def\epsilon{\varepsilon}

\newtheorem{theorem}{Theorem}[section]

\newtheorem{proposition}{Proposition}[section]
\newtheorem{lemma}{Lemma}[section]
\newtheorem{definition}{Definition}[section]
\newcommand{\cF}{{\cal F}}

\title{Sparsity Equivalence of Anisotropic Decompositions}
\author{Gitta Kutyniok}

\begin{document}
\maketitle

\begin{abstract}
Anisotropic decompositions using representation systems such as curvelets, contourlet, or shearlets
have recently attracted significantly increased attention due to the fact that they were shown to provide
optimally sparse approximations of functions exhibiting singularities on lower dimensional embedded manifolds.
The literature now contains various direct proofs of this fact and of related sparse approximation results.
However, it seems quite cumbersome to prove such a canon of results for each system separately, while many
of the systems exhibit certain similarities.

In this paper, with the introduction of the concept of {\em sparsity equivalence}, we aim to provide
a framework which allows categorization of the ability for sparse approximations of representation systems.
This framework, in particular, enables transferring results on sparse approximations from one system to another.
We demonstrate this concept for the example of curvelets and shearlets, and discuss how this viewpoint
immediately leads to novel results for both systems.
\end{abstract}

\vspace{.1in}
\vspace{.1in}

\vspace{.1in} {\bf Acknowledgements.} The author would like to thank Peter Binev, Emmanuel Cand\`{e}s,
Wolfgang Dahmen, Philipp Grohs, Demetrio Labate, Wang-Q Lim, and Pencho Petrushev for numerous discussions on this and related topics.
Special thanks go to David Donoho for enlightning comments and suggestions which helped to improve this
work. She would also like to thank the Department of Statistics at Stanford University and the Department
of Mathematics at Yale University for their hospitality and support during her
visits. This work was partially supported by Deutsche Forschungsgemeinschaft (DFG) Heisenberg fellowship KU
1446/8 as well as DFG Grants KU 1446/13 and KU 1446/14. \vspace{.1in}

\pagebreak


\section{Introduction}

Recently, a paradigm shift could be observed in applied mathematics, computer science, and electrical engineering.
The novel paradigm of sparse approximations now enables not only highly efficient encoding of functions
and signals, but also provides intriguing new methodologies, for instance, for recovery of missing data or
separation of morphologically distinct components.
At about the same time, scientists began to question whether wavelets are indeed perfectly suited for image
processing tasks, the main reason being that images are governed by edges while wavelets are isotropic objects.
This mismatch becomes also evident when recalling that Besov spaces can be characterized by the decay of
wavelet coefficient sequences however Besov models are clearly deficient to adequate capturing of edges.

These two fundamental observations have led to the research area of geometric multiscale analysis whose
main goal is to develop representation systems, preferably containing different scales, which are sensitive
to anisotropic features in functions/signals and provide sparse approximations of those. Such representation
systems shall for now be loosely coined {\em anisotropic systems}. Let us state as a few samples on the long list the directional
filter banks \cite{BS92}, directional wavelets \cite{AMV99}, ridgelets \cite{CD99b}, complex wavelets
\cite{Kin01}, (first and second generation) curvelets  \cite{CD04,CD05a,CD05b}, contourlets \cite{DV05},
bandlets \cite{MP05}, and shearlets \cite{GKL06,KL07}. Browsing through the literature, it becomes evident
that sparse approximation properties are quite similar for some systems such as curvelets and shearlets,
whereas other systems such as ridgelets show a different behavior. Delving more into the literature we observe that
for those systems exhibiting similar sparsity behavior many results were proven with quite resembling proofs.
One might ask: Is this cumbersome close repetition of proofs really necessary? We believe that the answer
is {\em no} and that a formalization of sparse approximation properties of anisotropic systems solves this
problem.

The main goal of this paper is to proclaim the concept of {\em sparsity equivalence} for anisotropic systems
leading to equivalence classes for sparsity properties, and thereby aiming for the aforementioned formalization
of sparse approximation properties. Our theoretical considerations are anticipated to have the following impacts:
\bitem
\item A thorough understanding of the ingredients of anisotropic systems which are crucial for
an observed sparse approximation property, thereby also categorizing different sparsity behaviors.
\item A framework within which sparsity results can be directly transferred from one system to others.
\item A quality measure for new anisotropic systems which they have to pass to be
considered eligible for a particular sparsity analysis.
\eitem


\subsection{The Concept of Sparsity Equivalence of Frame Expansions}
\label{subsec:concept}

Frame expansions are extensively utilized in applied mathematics, computer science, and electrical
engineering if non-uniqueness, yet stability is required, and might be regarded as a natural generalization
of the concept of an orthonormal basis. Non-uniqueness of an expansion is customarily exploited for deriving
resilience against erasures or quantization. However, lately the flexibility of such non-unique expansions has
been shown to lead to optimally sparse approximations of particular model classes of functions,
where sparsity of a coefficient sequence $(c_i)_{i \in I}$ is ideally measured in the $\| \, \cdot \, \|_0$-norm
counting the number of non-zero entries. The fundamental fact that this measure can be approximated by the
$\| \, \cdot \, \|_1$-norm as the closest convex norm has initiated and led to a deluge of results in the area
of sparse approximations and recovery; see the survey paper \cite{BDE09}.

Before continuing, let us briefly illustrate the precise relation of this sparsity measure with sparse approximation
properties. Given a tight frame $(\varphi_i)_{i \in I}$ for a Hilbert space $\cH$, say, and let $\cK \subset \cH$
be a class whose elements we desire to sparsely approximate. Approximation theory then paves the way to measure the
ability of  $(\varphi_i)_{i \in I}$ for sparse approximations of elements of $\cK$, and typically the decay of the squared error of
the `best' $n$-term approximation, i.e., the behavior of
\beq \label{eq:error}
\norm{f-\sum_{n \ge N} (\ip{f}{\varphi_i})_{(n)} \varphi_i}^2 \quad \mbox{as } N \to \infty,
\eeq
where $(\ip{f}{\varphi_i})_{(n)}$ is the $n$-th largest coefficient, is analyzed. Intriguingly, in the case of a
redundant system, it is not clear whether this is indeed the {\em best} $n$-term approximation; nevertheless it is
customarily exploited as a suitable substitute in lack of a more accurate and still conveniently applicable
selection rule. The term in \eqref{eq:error} can now be estimated by
\beq \label{eq:estimate}
\norm{f-\sum_{n \ge N} (\ip{f}{\varphi_i})_{(n)} \varphi_i}^2 \le C \cdot \sum_{n \ge N} |(\ip{f}{\varphi_i})_{(n)}|^2.
\eeq
Then the relation to $\|(\ip{f}{\varphi_i})_i\|_p$ ($0 < p \le 1$) is established by observing that
$\|(\ip{f}{\varphi_i})_i\|_p \le C'$ implies that the number of coefficients $\ip{f}{\varphi_i}$ exceeding $1/n$ is
bounded by $C' n^{1/p}$, thus the magnitude of the $n$-th largest coefficient $(\ip{f}{\varphi_i})_{(n)}$ is not
bigger than $C'' n^{-1/p}$.

As we already elaborated upon before, there do exist frames which show very similar sparse approximation
properties. Aiming towards a categorization of sparsity properties, we immediately observe that the well-exploited unitary
equivalence of frames does not serve our purposes here; the reason being that $\|(\ip{f}{U\varphi_i})_i\|_p =
\|(\ip{U^{-1} f}{\varphi_i})_i\|_p$ for all $f \in \cK$, however the class $\cK$ does not need to be invariant under
the unitary operator $U^{-1}$. Evidently, the equivalence relation we truly aim for is as follows:

\begin{definition}
Let $(\varphi_i)_{i \in I}$ and $(\psi_j)_{j \in J}$ be two frames for a Hilbert space $\cH$, let $\cK$
be a subset of $\cH$, and let $0 < p \le 1$. Then $(\varphi_i)_{i \in I}$ and $(\psi_j)_{j \in J}$ are {\em sparsity equivalent in $\ell_p$
with respect to $\cK$}, if, for each $f \in \cK$, we have $\|(\ip{f}{\varphi_i})_i\|_p < \infty$ if and only if
$\|(\ip{f}{\psi_j})_j\|_p < \infty$.
\end{definition}

This property is in fact a property of the cross-Grammian matrix $(\ip{\varphi_i}{\psi_j})_{i,j}$, more precisely, of
diagonal dominance of this matrix. A suitable norm for measuring the decay of this matrix away from the diagonal
was introduced in \cite{CD05b}, and is defined as follows: For $p \in (0,1]$, the $\| \cdot \|_{Op,p}$-norm of a matrix
$M=(m_{i,j})_{i,j}$ is given by
\[
\|M\|_{Op,p} = \max \Bigg\{ \Big(\sup_i\sum_j|m_{i,j}|^p\Big)^{1/p}, \Big(\sup_j\sum_i|m_{i,j}|^p\Big)^{1/p}\Bigg\}.
\]

This norm indeed measures whether sparsity equivalence is present, and we obtain the following result. Notice however,
that the condition on the cross-Grammian matrix is by far not necessary, which can be seen by the fact that it implies
sparsity equivalent in $\ell_p$ with respect to {\em any} subset $\cK$.

\begin{lemma} \label{lem:pnormsparsity}
Let $(\varphi_i)_{i \in I}$ and $(\psi_j)_{j \in J}$ be two tight frames for a Hilbert space $\cH$, let $\cK$ be a subset of $\cH$,
and let $0 < p \le 1$. If $\|(\ip{\varphi_i}{\psi_j})_{i,j}\|_{Op,p}$ is finite, then $(\varphi_i)_{i \in I}$ and
$(\psi_j)_{j \in J}$ are sparsity equivalent in $\ell_p$ with respect to $\cK$.
\end{lemma}

\noindent
{\bf Proof.}
Let $f \in \cC$, and assume that $\norm{(\ip{f}{\varphi_i})_i}_p < \infty$. From  $\|(\ip{\varphi_i}{\psi_j})_{i,j}\|_{Op,p}< \infty$
it follows that
\beq \label{eq:success}
\sup_i \sum_j |\ip{\varphi_i}{\psi_j}|^p < \infty \quad \mbox{and} \quad \sup_j \sum_i |\ip{\varphi_i}{\psi_j}|^p < \infty.
\eeq
Using the fact that $(\varphi_i)_i$ is a tight frame,
\[
\norm{(\ip{f}{\psi_j})_j}_p^p = \norm{(\langle\sum_i \ip{f}{\varphi_i}\varphi_i,\psi_j\rangle)_j}_p^p
= \norm{(\sum_i \ip{f}{\varphi_i}\langle\varphi_i,\psi_j\rangle)_j}_p^p.
\]
Now, since $p \le 1$,
\[
\norm{(\sum_i \ip{f}{\varphi_i}\langle\varphi_i,\psi_j\rangle)_j}_p^p
\le \sum_j \sum_i |\langle f,\varphi_i \rangle|^p \cdot |\langle\varphi_i,\psi_j\rangle|^p
\le \sum_i |\langle f,\varphi_i \rangle|^p \cdot \sup_i \sum_j |\langle\varphi_i,\psi_j\rangle|^p,
\]
which is finite by \eqref{eq:success} and due to the fact that $\norm{(\ip{f}{\varphi_i})_i}_p < \infty$.

For symmetry reasons, the implication $\norm{(\ip{f}{\psi_j})_j}_p < \infty \Rightarrow \norm{(\ip{f}{\varphi_i})_i}_p < \infty$
can be derived similarly. The lemma is proved.
\qed

We will now demonstrate this concept for the pair of curvelets and shearlets, which are two prominent examples of
anisotropic systems even sharing parabolic scaling as the main anisotropic force. The intuition that they
should be sparsity equivalent is substantiated by comparing results on sparse approximation properties of curvelets
and shearlets. And, in fact, the result derived in Subsection \ref{subsec:main} shows this to be true. Before stating
the result, we first need to introduce those two systems.


\subsection{Curvelets and Shearlets}
\label{subsec:curveshear}

We now recall the definitions of curvelets -- focussing on second generation curvelets -- and shearlets. Those two
systems will be exemplarily focused on in our demonstration of the framework of sparsity equivalence.


\subsubsection{Curvelets}
\label{subsec:2gencurve}

The main motivation for the introduction of curvelets came from the observation that -- by taking a computer
vision point of view -- edges are those features governing an image while separating smooth regions. A
first model for this view point was introduced in \cite{Don99} and coined a `cartoon-like model'. This
model then in fact revealed the suboptimal treatment of edges by the at that time seemingly superior system of
wavelets.

The introduction of (first generation) tight curvelet frames in 2004 by Cand\'{e}s and Donoho \cite{CD04}, which
provably provide (almost) optimally sparse approximations within such a cartoon-like model might be considered a
milestone in applied harmonic analysis. Later, second generation curvelets were introduced in \cite{CD05b} due to a
more satisfactory associated system with continuous parameters \cite{CD05a}, and were shown to provide optimally
sparse decompositions of Fourier Integral Operators \cite{CD05}.

To present the definition of these second generation curvelets -- from now on also called {\em curvelets} in
contrast to {\em first generation curvelets} --, let $W$ be the Fourier transform of a one-di\-men\-sio\-nal wavelet and $V$ be a
`bump function' in Fourier space. We select both functions to be band-limited, where $\Sp W \subseteq [-2,-1/2] \cup [1/2,2]$ and
$\Sp V \subseteq [-1,1]$, and to satisfy $W, V \in C^\infty$. Curvelets live on anisotropic regions of width $2^{-j}$
and length $2^{-j/2}$ at various orientations, which are parameterized by angle. For our purposes, it is sufficient to
ignore the low frequency part in curvelet decompositions as discussed latter. We just mention that appropriate
low frequency functions can be added to the  curvelet system defined below to force it to become a tight frame for
$L^2(\RR^2)$. Hence we will only state the definition of curvelets restricted to
\[
\cC = \{\xi \in \RR^2 : \norm{\xi}_\infty \ge 1\}.
\]

\begin{figure}[ht]
\begin{center}
\includegraphics[height=1.75in]{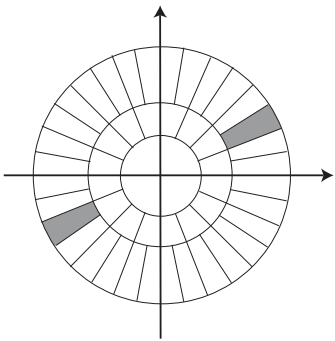}
\end{center}
\caption{The tiling of the frequency domain induced by curvelets.}
\label{fig:Curvelets}
\end{figure}

Let now $A_a$ denote the parabolic scaling matrix $A_a = \diag(a, \sqrt{a})$. {\em Curvelets} at scale $j \ge 0$, orientation
$\ell = 0, \dots, 2^{j/2}-1$, and spatial position $m = (m_1,m_2) \in \bZ^2$ are then defined by their Fourier transforms
of some $\xi \in \RR^2$, with $(r,\omega)$ denoting the associated polar coordinates,
\[
      \hat{\gamma}_{\mu}(\xi) =  2^{-j\frac{3}{4}}  \cdot W(r/2^{j}) V((\omega-\theta_{j,\ell})2^{j/2})
         \cdot e^{ i \langle R_{\theta_{j,\ell}}A_{2^{-j}}m,\xi\rangle },
\]
where here $\theta_{j,\ell} = 2\pi \ell /2^{j/2}$, $R_\theta$ is planar rotation by $- \theta$ radians, and we let $\mu = (j,\ell,m)$
index scale, orientation, and position. We refer to \cite[Sect. 4.3, pp. 210-211]{CD05b} for more details,
and to Figure \ref{fig:Curvelets} for an illustration of the induced tiling of the frequency plane.

\subsubsection{Shearlets}
\label{subsec:shearlets}

In 2006, a novel directional representation system -- so-called shearlets -- has been proposed in
\cite{GKL06,KL07}, which provides a unified treatment for the continuum and digital world. The main point in
comparison with curvelets is the fact that angles are replaced by slopes when parameterizing directions which greatly
supports the treating of the digital setting. Hence the theory of shearlets allows an associated digital theory which
can be directly implemented \cite{DKS08}.

\begin{figure}[ht]
\begin{center}
\includegraphics[height=1.75in]{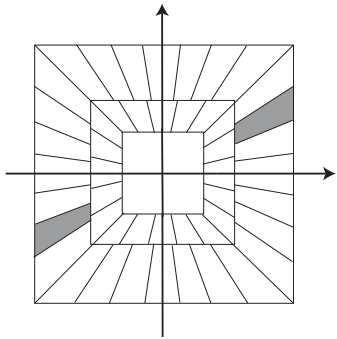}
\end{center}
\caption{The tiling of the frequency domain induced by cone-adapted shearlets.}
\label{fig:ShearletsCone}
\end{figure}

In a similar way as curvelets do shearlets live on anisotropic regions of width $2^{-j}$ and length $2^{-j/2}$ at various
orientations, which are now parameterized by slope rather than angle as for curvelets. Similar to the definition of
curvelets stated in Subsection \ref{subsec:2gencurve}, also here we will ignore the low frequency part, and just mention that it can be appropriately included to
yield a tight frame for $L^2(\RR^2)$. Let now the Fourier transform $W$ of a wavelet and a bump function $V$ be chosen
as in Subsection \ref{subsec:2gencurve}, and let $\cC^{(1)}$ and $\cC^{(2)}$ denote the following two cones:
\[
\cC^{(\iota)} = \left\{ \begin{array}{rcl}
\{(\xi_1,\xi_2) \in \bR^2 : |\xi_1| \ge 1,\, |\xi_2/\xi_1| \le 1\} & : & \iota = 1,\\
\{(\xi_1,\xi_2) \in \bR^2 : |\xi_2| \ge 1,\, |\xi_1/\xi_2| \le 1\} & : & \iota = 2.
\end{array}
\right.
\]
For cone $\cC^{(1)}$, at scale $j \ge 0$, orientation $k = -\lceil 2^{j/2}\rceil, \dots,$ $\lceil 2^{j/2}\rceil$, and spatial position $m \in \bZ^2$,
the associated {\em shearlets} are defined by their Fourier transforms
\[
\hat{\sigma}_{\eta}(\xi) = 2^{3j/4} {\psi}({S}_{k} {A}_{2^j}\cdot-m),
\]
where $S_k$ denotes the shear matrix
\[
S_k = \begin{pmatrix}
  1 & k \\ 0 & 1
\end{pmatrix},
\]
and $\eta = (j,k,m,1)$ indexes scale, orientation, position, and cone. We now assume that $\psi \in L^2(\RR^2)$
is chosen such that
\[
\hat{\psi}(\xi_1,\xi_2) = W(\xi_1) V(\xi_2/\xi_1),
\]
wherefore
\[
\hat{\sigma}_{\eta}(\xi) = 2^{-j\frac{3}{4}} W(\xi_1/2^j) V(k + 2^{j/2}\xi_2/\xi_1) e^{i\langle S_k^T A_{2^{-j}}m,\xi\rangle}.
\]
The shearlets for $\cC^{(2)}$ are defined likewise by symmetry, as illustrated in Figure \ref{fig:ShearletsCone};
this initiated the terminology {\em cone-adapted shearlets} in contrast to shearlets arising directly from a group
representation (cf. \cite{KKL10b}).

We remark that the discrete shearlets considered, for instance, in \cite{GL07a} differ slightly from this choice, since they are usually
associated with a scaling of $4^j$.
However, it is easily checked -- and we refer concerning this issue and additional details to the survey paper \cite {KLL10b}~-- that
the shearlets as defined here also form a tight frame for $L^2(\RR^2)$.

The attentive reader will have also noticed that we here consider the class of band-limited shearlets although there
has just recently been introduced a class of compactly supported shearlets which have superior spatial domain
localization (see \cite{KL10,KKL10a}). Since in this paper we however aim to compare curvelets and shearlets and since curvelets
are band-limited, the class of band-limited shearlets is the canonical choice. Another issue to consider is the
fact that compactly supported shearlets are not a tight frame, thereby requiring adaptions to the analysis.
Additional thoughts on compactly supported versus band-limited shearlets can be found in Section \ref{sec:extensions}.


\subsection{Equivalence Result}
\label{subsec:main}

The introduction of the concept of sparsity equivalence in Subsection \ref{subsec:concept} now motivates us
to ask whether curvelets and shearlets belong to the same equivalence class, hence are sparsity equivalent.
The many quite similar results on sparse approximation properties of those two systems seem to indicate this.
According to Lemma \ref{lem:pnormsparsity}, the $\ell_p$ norm of the cross-Grammian matrix reveals the
true sparsity relation, and we obtain the following result, whose lengthy proof is presented in
Subsection \ref{subsec:proof}.

\begin{theorem} \label{theo:main1}
For all $0 < p \le 1$,
\[
\norm{(\ip{\sigma_{\eta}}{\gamma_{\mu}})_{\eta, \mu}}_{Op,p} < \infty.
\]
\end{theorem}


Now Lemma \ref{lem:pnormsparsity} can be applied to derive the already intuitively expected sparsity equivalence
of shearlets and curvelets.

\begin{theorem} \label{theo:main2}
For all $0 < p \le 1$, the shearlet frame $(\sigma_\eta)_\eta$ and the curvelet frame $(\gamma_\mu)_\mu$ are
sparsity equivalent in $\ell_p$ with respect to $L^2(\RR^2)$.
\end{theorem}

%


\subsection{Impact of Sparsity Equivalence}

The significance of the viewpoint of sparsity equivalence lies in the fact that it not only provides a
thorough understanding of the ability of different anisotropic systems for sparse expansions when compared
to each other -- thereby providing a qualitative comparison --, but it moreover allows the transfer
of sparsity results without repeating quite similar proofs.

The theorem presented in the previous subsection is a first demonstration of the power of such a higher
level viewpoint of sparsity behavior. In fact, this result automatically leads to novel results on and insights
in sparse expansions by curvelets and shearlets. A few examples, for which this conceptually
new approach is fruitful, will be presented in Section \ref{sec:impact} including optimally sparse
approximations of cartoon-like images and the ability for geometric separation
of morphologically distinct phenomena.


\subsection{Extensions and General Viewpoint}
\label{subsec:extensions}

As mentioned before, Theorems \ref{theo:main1} and \ref{theo:main2} are amenable to generalizations and extensions.
Previewing Section \ref{sec:extensions}, we briefly discuss a few examples.

\bitem
\item {\em Curvelets and Shearlets}. A similar statement as Theorem \ref{theo:main2} should be provable for
first generation curvelets as also for the new class of compactly supported shearlets.
\item {\em Other Systems}. The analysis of sparsity equivalence of curvelets and shearlets we drove here can
and should be applied to other pairs of systems. Ideally, novelly introduced systems could be compared to
a system whose sparse approximation properties are already very well understood.
\item {\em Systems with Continuous Parameters}. Certainly, we can also ask about similar sparsity
behavior for systems with continuous parameters. This however requires a different sparsity model; one
conceivable path would be to compare their ability to resolve wavefront sets.
\item {\em Weighted Norms}. When aiming at transferring results such as sparse decompositions of curvilinear
integrals \cite{CD00a} or sparse decompositions of the Radon transform \cite{CD00b}, sometimes weighted
$\ell_p$ norms might need to be analyzed. This is also essential for analyzing associated approximation
spaces.
\eitem



\subsection{Outline}

We start by presenting the analysis of sparsity equivalence between curvelets and shearlets and providing
the proof of Theorem \ref{theo:main1}. We then analyze the impact of this and related results on sparse
approximation properties of anisotropic systems in Section \ref{sec:impact}.  In particular, we derive
novel results on sparse approximation of cartoon-like images using curvelets and on the ability of geometric
separation using shearlets and wavelets. This section is followed by a discussion on extensions of our
framework (see Section \ref{sec:extensions}).


\section{Sparsity Equivalence between Curve\-lets and Shearlets}
\label{sec:curveshear}

In this section our goal is to prove sparsity equivalence in $\ell_p$ of curvelets and
shearlets for all $0 < p \le 1$. Due to Lemma \ref{lem:pnormsparsity}, this task is reduced to proving Theorem \ref{theo:main1}, i.e.,
showing that the $\| \cdot \|_{Op,p}$-norm of the cross-Grammian matrix of curvelets and shearlets is finite.

We first realize that for our analysis we only need to consider those curvelets and shearlets
which respond to the high-frequency content of a function. More precisely, if we are given
a function, say $f \in L^2(\RR^2)$, we might decompose it as $f = f_L + f_H = g_L \cdot f + g_H \cdot f$,
where $g_L$ is a low pass filter with $(\hat{g}_L)|_\cC \equiv 1$, and $g_H$ is an `associated' high pass
filter satisfying $g_L + g_H = 1$. Now notice, that the inner products between elements of both frames
corresponding to $g_L$ are negligible due to their almost orthogonality, since they are scaling functions;
also the inner products of those elements with elements corresponding to $g_H$ are of a similar reason
negligible.

This argument shows that it is sufficient to only consider the cross-Grammian matrix of the elements of the
curvelet and shearlet frame introduced in Subsection \ref{subsec:curveshear}, i.e., those analyzing
the high-frequency part of a function.


\subsection{Estimates for the Entries of the Cross-Grammian Matrix}

We start by establishing estimates on the absolute values of inner products of curvelets and
shearlets. An essential ingredient will be the following well-known result, which we state
here for the convenience of the reader. A detailed proof might for instance be found in
\cite[Lem. 2.3]{KL07}.

\begin{lemma}
Suppose $g$ satisfies $\hat g \in C_0^\infty(\RR^d)$ with $\hat g$ being supported
on a fixed bounded rectangle $R \subset \R^d$. Then, for each $N \in \N$, there exists a
constant $C_N$ such that
\[
|g(x)| \le C_N \, (1+|x|^2)^{-N}\quad \mbox{for all } x \in \RR^d.
\]
In particular, $C_N = N \, \lambda(R) \, \bigl(\norm{\hat g}_\infty + \norm{\Delta^N \hat g}_\infty \bigr)$,
where $\Delta = \sum_{i=1}^d \frac{\partial^2}{\partial \xi_i^2}$ denotes the frequency
domain Laplacian operator and $\lambda(R)$ is the Lebesgue measure of $R$.
\end{lemma}

In \cite{CD05b} the following conclusion was drawn from this lemma which we will also require for our proof.

\begin{lemma}\cite[Lem. 5.6]{CD05b} \label{lem:decaybandlimited}
Suppose $(f_j)_{j \ge 0}$ is a sequence of functions satisfying that each $\hat{f}_j$ is supported
in a rectangle $R_j = A_{2^{j}}([-C_1,C_1] \times [-C_2,C_2])$ and every scaled function
\[ \hat{g}_j(\xi) = 2^{\frac32 j} \hat{f}_j(A_{2^{j}}\xi)\]
obeys $\norm{\hat{g}_j}_{C^N} \le \rho_N$ for $N=2,4,6,\ldots$ with each $\rho_N$ being independent
on $j$. Then, for $N=2,4,6,\ldots$, there exist constants $C_N$ such that
\[
|f_j(x)| \le C_N (\rho_0 + \rho_N)\langle A_{2^j} x\rangle^{-N}\quad \mbox{for all } x \in \RR^d,
\]
where
\[
\langle y \rangle = (1+y^2)^{1/2}.
\]
\end{lemma}

The estimates which are proved in the following proposition are carefully designed so that
the previously stated claim concerning the $\| \cdot \|_{Op,p}$-norm, $0 < p \le 1$ of the cross-Grammian
matrix of curvelets and shearlets does follow almost immediately as a corollary. We note that a similar
estimate for the second cone $\cC^{(2)}$ holds with a resembling proof.

\begin{proposition} \label{prop:CrossMatrixMainEstimate}
Let $j,\t{j} \ge 0$, $|k| \le \lceil2^{j/2}\rceil$, $0 \le \ell < 2^{j/2}$, and $m,\t{m} \in \ZZ$. Then, for each
$N=2,4,6,\ldots$, there exist constants $C_N$ so that
\[
\absip{\sigma_{j,k,m,1}}{\gamma_{\t{j},\ell,\t{m}}} \le C_N 1_{\{|j-\t{j}|\le 2\}}
1_{\{k \in K_{j,\t{j},\ell}\}} 1_{\{\ell \in L_{j,\t{j},k}\}} \langle |b_{j,k,m,\t{j},\ell,\t{m}}| \rangle^{-N}
\]
where
\[
K_{j,\t{j},\ell} = \{k : \lfloor -2^{j/2}\cdot\tan(2^{-\t{j}/2}(1+2\pi\ell))-1\rfloor \le k \le \lceil -2^{j/2}\cdot\tan(2^{-\t{j}/2}(-1+2\pi\ell))+1\rceil\},
\]
\[
L_{j,\t{j},k} = \{\ell :  \lfloor2^{\t{j}/2}\arctan(2^{-j/2}(-1-k))-1\rfloor \le 2\pi\ell \le \lceil 2^{\t{j}/2}\arctan(2^{-j/2}(1-k))+1\rceil\},
\]
and
\[
b_{j,k,m,\t{j},\ell,\t{m}} = A_{2^{j}}(S_k^TA_{2^{-j}}m-R_{\theta_{\t{j},\ell}}A_{2^{-\t{j}}}\t{m}).
\]

\end{proposition}

\begin{figure}[ht]
\centering
\includegraphics[height=2.3in]{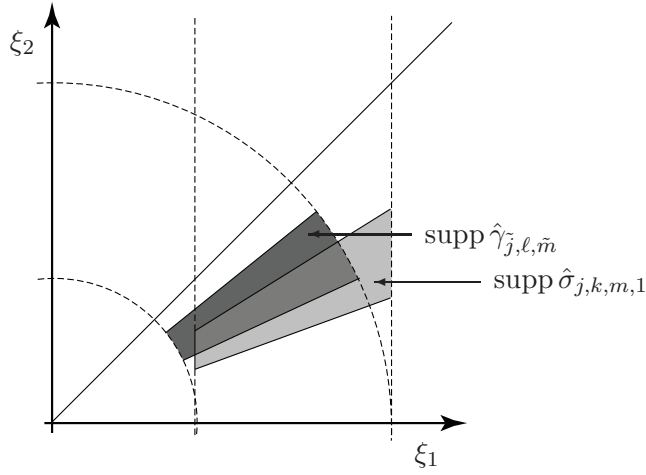}
\put(5,60){\vector(-1,0){40}}
\put(10,58){$\Sp \hat \sigma_{j,k,m,1}$}
\put(-21,78){\vector(-1,0){39}}
\put(-16,75){$\Sp  \hat \gamma_{\t{j},\ell,\t{m}}$}
\put(-20,-8){$\xi_1$}
\put(-173,148){$\xi_2$}
\caption{The support in frequency domain of some functions $\sigma_{j,k,m,1}$ and
$\gamma_{\t{j},\ell,\t{m}}$.}
\label{fig:1}
\end{figure}

\noindent
{\bf Proof.}
To illustrate the different supports in frequency domain of $\sigma_{j,k,m,1}$ and $\gamma_{\t{j},\ell,\t{m}}$,
a property which will be exploited in the sequel, we refer to Figure \ref{fig:1}.

We now fix $j,k,m$. By employing Plancherel's theorem, we have
\beq \label{eq:FourierEstimate1}
\absip{\sigma_{j,k,m,1}}{\gamma_{\t{j},\ell,\t{m}}}
= 2^{-(j+\t{j})\frac{3}{4}}  \int f_{j,\t{j},k,\ell}(\xi) \cdot
e^{ i \langle S_k^TA_{2^{-j}}m-R_{\theta_{\t{j},\ell}}A_{2^{-\t{j}}}\t{m},\xi\rangle }\, d\xi,
\eeq
where
\[
 \hat f_{j,\t{j},k,\ell}(\xi) = W(\xi_1/2^j)W(r/2^{\t{j}})  V(k + 2^{j/2}\xi_2/\xi_1)
        V((\omega-\theta_{\t{j},\ell})2^{\t{j}/2}).
\]
Due to the support conditions of $W$ and $V$, the support of $\hat \sigma_{j,k,m,1}$
equals
\beq \label{eq:SupportShearlet1}
\Sp \hat \sigma_{j,k,m,1} = \{(\xi_1,\xi_2) \in \RR^2 : \xi_1 \in [2^{j-1},2^{j+1}], \,
\xi_2/\xi_1 \in 2^{-{j/2}}([-1,1]-k)\},
\eeq
whereas the support of $\hat \gamma_{\t{j},\ell,\t{m}}$ is
\beq \label{eq:SupportCurvelet1}
\Sp \hat \gamma_{\t{j},\ell,\t{m}} = \{(\xi_1,\xi_2) \in \RR^2 : r \in [2^{\t{j}-1},2^{\t{j}+1}],\,
\omega \in 2^{-{\t{j}/2}}[-1,1]+\theta_{\t{j},\ell}\}.
\eeq
We conclude that $\hat f_{j,\t{j},k,\ell} \equiv 0$ unless $|j-\t{j}|\le 2$, hence
\beq \label{eq:EstimateCross1}
\absip{\sigma_{j,k,m,1}}{\gamma_{\t{j},\ell,\t{m}}} \le C_{j,k,m,\t{j},\ell,\t{m}} 1_{\{|j-\t{j}|\le 2\}}.
\eeq
Our next task is to estimate the range of $\ell$ for which $\absip{\sigma_{j,k,m,1}}{\gamma_{\t{j},\ell,\t{m}}}$
is non-zero. This will
be done by showing that this parameter is contained in a compact set whose size
is uniformly bounded as $j, \t{j} \to \infty$. For this, we will study the slopes of the boundaries of the
supports of $\hat \sigma_{j,k,m,1}$ and $\hat \gamma_{\t{j},\ell,\t{m}}$ in angular direction.
For better comparison with \eqref{eq:SupportShearlet1}, the support \eqref{eq:SupportCurvelet1} might be rewritten as
\beq \label{eq:SupportCurvelet2}
\Sp \hat \gamma_{\t{j},\ell,\t{m}} = \{(\xi_1,\xi_2) \in \RR^2 : r \in [2^{\t{j}-1},2^{\t{j}+1}],\,
\xi_2/\xi_1 \in \tan(2^{-{\t{j}/2}}([-1,1]+2\pi \ell))\}.
\eeq
Notice that the angle between the two angular boundary lines of the support of curvelets does not change
with $\ell$, whereas in the shearlet case the angle becomes smaller as the support of the Fourier transform of
the  shearlet approaches the angle bisector of the first quadrant.
From \eqref{eq:SupportShearlet1}
and \eqref{eq:SupportCurvelet2}, it follows that $\hat f_{j,\t{j},k,\ell} \equiv 0$
if
\[
\tan(2^{-{\t{j}/2}}(1+2\pi \ell)) \le 2^{-{j/2}}(-1-k)\quad \mbox{or} \quad
\tan(2^{-{\t{j}/2}}(-1+2\pi \ell)) \ge 2^{-{j/2}}(1-k).
\]
Continuing \eqref{eq:EstimateCross1}, this implies
\beq \label{eq:EstimateCross2}
\absip{\sigma_{j,k,m,1}}{\gamma_{\t{j},\ell,\t{m}}} \le C_{j,k,m,\t{j},\ell,\t{m}} 1_{\{k \in K_{j,\t{j},\ell}\}}
1_{\{\ell \in L_{j,\t{j},k}\}},
\eeq
with $K_{j,\t{j},\ell}$ and $L_{j,\t{j},k}$ as defined in the statement of the lemma.

Next we aim to estimate the decay in $m$ and $\t{m}$ by making use of Lemma \ref{lem:decaybandlimited}.
To prepare the application of this lemma, we rescale the function $f_{j,\t{j},k,\ell}$ in the
term of the RHS of \eqref{eq:FourierEstimate1} according to
\[
\hat g_{j,\t{j},k,\ell}(u,v) = 2^{j \frac32}\hat{f}_{j,\t{j},k,\ell}(A_{2^j}(u,v)).
\]
This yields a function which can be decomposed into factors in the following way:
\[
\hat g_{j,\t{j},k,\ell}(u,v) = \tilde{W}_{0,j}(u,v)\tilde{W}_{1,j}(u,v)\tilde{V}_{0,j}(u,v)\tilde{V}_{1,j}(u,v).
\]
All factors belong to $C^\infty$, and it can be checked that their derivatives are
bounded independent on $j$ (for a similar argument confirm \cite[Subsec. 5.2]{CD05b}).
This allows us to apply Lemma \ref{lem:decaybandlimited} to obtain
\[
|f_{j,\t{j},k,\ell}(b)| \le c_N \langle |A_{2^{j}}b| \rangle^{-N}  ,\qquad N=2,4,6,\ldots.
\]
From this we conclude that, for $N=2,4,6,\ldots$,
\beqn
|\langle \sigma_{j,k,m,1},\gamma_{\t{j},\ell,\t{m}}\rangle|
& = & |f_{j,\t{j},k,\ell}(S_k^TA_{2^{-j}}m-R_{\theta_{\t{j},\ell}}A_{2^{-\t{j}}}\t{m})|\\
& \le & c_N \langle |A_{2^{j}}(S_k^TA_{2^{-j}}m-R_{\theta_{\t{j},\ell}}A_{2^{-\t{j}}}\t{m})| \rangle^{-N}.
\eeqn
Combining this estimate with the estimates from \eqref{eq:EstimateCross1} and \eqref{eq:EstimateCross2}
proves the lemma.
\qed


\subsection{Proof of Theorem \ref{theo:main1}}
\label{subsec:proof}

Let $0 < p \le 1$. We start by proving that
\beq \label{eq:pNormFirstPart}
\sup_\mu \sum_\eta |(\ip{\sigma_{\eta}}{\gamma_{\mu}})_{\eta, \mu}|^p < \infty.
\eeq
Setting $\eta = (j,k,m,1)$  and $\mu = (\t{j},\ell,\t{m})$, by Proposition
\ref{prop:CrossMatrixMainEstimate},
\begin{eqnarray} \nonumber
\sup_\mu \sum_\eta |(\ip{\sigma_{\eta}}{\gamma_{\mu}})_{\eta, \mu}|^p
& \le & C_{N,p} \sup_{\t{j},\ell,\t{m}} \sum_{\{|j-\t{j}|\le 2\}} \sum_{\{k \in K_{j,\t{j},\ell}\}}
\sum_m \langle |b_{j,k,m,\t{j},\ell,\t{m}}| \rangle^{-pN}\\ \label{eq:EstimatePNorm1}
& \le & C_{N,p}' \sup_{j,\ell,\t{m}} \sum_{\{k \in K_{j,j,\ell}\}}
\sum_m \langle |b_{j,k,m,j,\ell,\t{m}}| \rangle^{-pN}.
\end{eqnarray}
The last estimate was derived by observing that the maximum of $\sum_m \langle |b_{j,k,m,j,\ell,\t{m}}| \rangle^{-pN}$ is
attained if $\tilde{j} = j$.

We next compute the number of integers $k$ satisfying $|k| \le \lceil2^{j/2}\rceil$ which are contained in $K_{j,j,\ell}$. We observe
that $\#(K_{j,j,\ell})$ is maximal if $\ell$ is chosen so that the upper bound of the
curvelet coincides with the angle bisector,
the reason being that the support in frequency domain of this `corner curvelet' has a maximal
number of intersections with frequency supports of shearlets. In fact, the angular support
of the Fourier Transform of shearlets
become smaller when the angle increases, hence more shearlets are needed to overlap
the angular frequency support of a curvelet, which does not change its size with varying angle
(also compare the proof of Proposition \ref{prop:CrossMatrixMainEstimate}).
Hence, using
\[
\Sp \hat \gamma_{\t{j},\ell,\t{m}} = \{(\xi_1,\xi_2) \in \RR^2 : r \in [2^{\t{j}-1},2^{\t{j}+1}],\,
\omega \in 2^{-{\t{j}/2}}[-1,1]+\theta_{\t{j},\ell}\}
\]
(cf. \eqref{eq:SupportCurvelet1}), it is sufficient to restrict to the situation
\[ \frac{\pi}{4} = 2^{-{j/2}}+\theta_{j,\ell}.\]
By definition of $\theta_{j,\ell}$, we therefore obtain the condition
\[
\ell=(2\pi)^{-1}(2^{j/2} \pi/4-1).
\]
The definition of $K_{j,j,\ell}$ implies
\beqn
\#(K_{j,j,\ell})
& \hspace*{-0.2cm} \le \hspace*{-0.2cm} & -2^{j/2} \tan(2^{-{j/2}}(-1+2^{j/2}\pi/4-1)) - (-2^{j/2} \tan(2^{-{j/2}}(1+2^{j/2}\pi/4-1))) + 3\\
& \hspace*{-0.2cm} \le \hspace*{-0.2cm} & 2^{j/2} (1-\tan(\pi/4-2\cdot 2^{-{j/2}})) + 3.
\eeqn
Now
\[
2^{j/2} (1-\tan(\pi/4-2\cdot 2^{-{j/2}})) \to 4, \qquad j \to \infty,
\]
hence,
\[
\#(K_{j,j,\ell}) \to 7, \qquad j \to \infty.
\]
From \eqref{eq:EstimatePNorm1}, we can then conclude that
\beq \label{eq:EstimatePNorm2}
\sup_\mu \sum_\eta |(\ip{\sigma_{\eta}}{\gamma_{\mu}})_{\eta, \mu}|^p
\le C_{N,p}'' \sup_{j,k,\ell,\t{m}}
\sum_m \langle |b_{j,k,m,j,\ell,\t{m}}| \rangle^{-pN}.
\eeq
Next we aim to prove that
\beq \label{eq:EstimateSumB}
\sum_m \langle |b_{j,k,m,j,\ell,\t{m}}| \rangle^{-pN}
\le C_{N,p} \sum_m \langle |m| \rangle^{-pN} \le C_{N,p}'.
\eeq
The second inequality follows easily from the facts that $\langle |m| \rangle^{-1}
\le 2^{-1} \langle m_1 \rangle^{-1} \langle m_2 \rangle^{-1}$ and choosing $pN$ large
enough such that $\sum_{m_i} \langle m_i \rangle^{-pN} < \infty$ for $i=1,2$.
Concerning the first inequality in \eqref{eq:EstimateSumB}, recall that
\[
b_{j,k,m,j,\ell,\t{m}} = A_{2^{j}}(S_k^TA_{2^{-j}}m-R_{\theta_{j,\ell}}A_{2^{-j}}\t{m}).
\]
Since we sum over $m$, WLOG we can assume that $\t{m}=0$. We have
\[
|A_{2^{j}}S_k^TA_{2^{-j}}m| = |(m_1,2^{-j/2}km_1+m_2)|,
\]
and hence it follows immediately that
\[
\sum_m \langle |A_{2^{j}}S_k^TA_{2^{-j}}m| \rangle^{-pN}
\le C_{N,p} \sum_m \langle |m| \rangle^{-pN}.
\]
This completes the proof of \eqref{eq:EstimateSumB}.

Finally, \eqref{eq:pNormFirstPart} follows from the application of \eqref{eq:EstimateSumB}
to \eqref{eq:EstimatePNorm2} and an estimate similar to Proposition \ref{prop:CrossMatrixMainEstimate}
for the second cone $\cC^{(2)}$ to handle the indices $\eta = (j,k,m,2)$.

It remains to prove that
\beq \label{eq:pNormSecondPart}
\sup_\eta \sum_\mu |(\ip{\sigma_{\eta}}{\gamma_{\mu}})_{\eta, \mu}|^p < \infty.
\eeq
Again, by Proposition \ref{prop:CrossMatrixMainEstimate},
\begin{eqnarray} \nonumber
\sup_\eta \sum_\mu |(\ip{\sigma_{\eta}}{\gamma_{\mu}})_{\eta, \mu}|^p
& \le & C_N^p \sup_{j,k,m} \sum_{\{|j-\t{j}|\le 2\}} \sum_{\{\ell \in L_{j,\t{j},k}\}}
\sum_{\t{m}} \langle |b_{j,k,m,\t{j},\ell,\t{m}}| \rangle^{-pN}\\ \label{eq:EstimatePNorm3}
& \le & C_N' \sup_{j,k,m} \sum_{\{\ell \in L_{j,\t{j},k}\}}
\sum_{\t{m}} \langle |b_{j,k,m,j,\ell,\t{m}}| \rangle^{-pN}.
\end{eqnarray}
We now need to estimate $\#(L_{j,j,k})$. Recalling our `worst-case-discussion'
in the previous case, the number of elements in $L_{j,j,k}$ reaches its maximum if $k=0$, i.e.,
the Fourier transform of the shearlet associated with $k$ `sits' precisely on
the $x$-axis. In this case, using the definition of $L_{j,j,k}$,
\beqn
\#(L_{j,j,k})
& \le & (2\pi)^{-1}(2^{j/2} \arctan(2^{-{j/2}}) - 1) - (2\pi)^{-1}(2^{j/2} \arctan(-2^{-{j/2}}) + 1)+3\\
& = & (2\pi)^{-1}2^{j/2}(\arctan(2^{-{j/2}})-\arctan(-2^{-{j/2}})) - \pi^{-1} + 3\\
& \le & \pi^{-1}(2^{j/2}\arctan(2^{-{j/2}}) - 1) + 3.
\eeqn
Consequently,
\[
\#(L_{j,j,k}) \to 3, \qquad j \to \infty.
\]
Hence, continuing the computation in \eqref{eq:EstimatePNorm3},
\[
\sup_\eta \sum_\mu |(\ip{\sigma_{\eta}}{\gamma_{\mu}})_{\eta, \mu}|^p
\le C_N'' \sup_{j,\ell,k,m}
\sum_{\t{m}} \langle |b_{j,k,m,j,\ell,\t{m}}| \rangle^{-pN}.
\]
Combining this estimate with
\[
\sum_{\t{m}} \langle |b_{j,k,m,j,\ell,\t{m}}| \rangle^{-pN}
\le C_{N,p} \sum_{\t{m}} \langle |\t{m}| \rangle^{-pN} \le C_{N,p}',
\]
which can be proven similarly as \eqref{eq:EstimateSumB} (cf. also \cite[Sect. 5.2]{CD05b}), the claim \eqref{eq:pNormSecondPart}
follows. This completes the proof.


\section{Impact of Sparsity Equivalence}
\label{sec:impact}

To illustrate the impact of the concept of sparsity equivalence focussing on the chosen exemplary
case of curvelets and shearlets, we now discuss two different situations in which the application
of Theorem \ref{theo:main2} automatically leads to novel results.

We might have also included the
search for optimally sparse expansions of Fourier Integral Operators of order 0. Since such a
result is however already known for  curvelets and shearlets -- with not surprisingly
quite similar proofs --, our considerations cannot lead to new results. They however point to a
simplified analysis once the result was known for either  curvelets or shearlets.


\subsection{Optimal Sparse Representation of $C^2$-Curvilinear Singularities}

To efficiently process image data, optimally sparse approximations are crucial. As already discussed in Subsection
\ref{subsec:concept}, the ability to sparsely approximate a class of signals is measured by the decay of the
error of the $n$-term approximation using the largest $n$ coefficients in magnitude; see \eqref{eq:error}.
Choosing the `correct' model class for images is certainly a highly delicate task. In 2004, Cand\`{e}s and
Donoho proposed a so-called cartoon model \cite{CD04} motivated by the fact that edges are the most prominent
features in images, a fact also evidenced in computer vision.

The cartoon model they proclaimed is defined
as follow: Let $B \subset [0,1]^2$ be bounded by a closed $C^2$ curve whose curvature is uniformly bounded
by some $\nu >0$, and let $STAR^2(\nu)$ be the class of translates of such sets $B$. Then the class of
{\em cartoon-like images} $\cE^2(\nu)$ is defined to be the set of functions $f$ on $\RR^2$ of the form
\[
f = f_0 + f_1 \chi_{B},
\]
where $f_0,f_1 \in C^2(\RR^2)$ with compact support in $[0,1]^2$, $B \in STAR^2(\nu)$, and $\|f\|_{C^2} =
\sum_{|\alpha| \leq 2} \|D^{\alpha}f\|_{\infty} \leq 1.$

By information theoretic arguments, it can be shown that the optimally achievable rate of sparse approximations
under weak conditions on the dictionary and the selection process is $N^{-2}$ as $N \to \infty$.
For first generation curvelets \cite{CD04} as well as for shearlets \cite{GL07a} (see also \cite{KL10}), this
rate is achieved up to a multiplicative log factor of $(\log N)^{3}$.

We now claim that also (second generation) curvelets achieve the optimal sparse approximation rate up to
a factor negligible compared to $N^{-2}$.

\begin{theorem}
The curvelet frame $(\gamma_\mu)_\mu$ provides (almost) optimally sparse approximations of functions
$f \in \cE^2(\nu)$, i.e., there exists some $C > 0$ such that
\[
\|f-f_N\|_2^2 \leq C \cdot N^{-2} \cdot D(N) \qquad \text{as } N \rightarrow \infty,
\]
where $f_N$ is the nonlinear N-term approximation obtained by choosing the N largest curvelet coefficients
of $f \in \cE^2(\nu)$ and $N^{-\epsilon} \cdot D(N) \to 0$ as $N \to 0$ for all $\epsilon > 0$.
\end{theorem}

\noindent
{\bf Proof.}
Given some $f \in \cE^2(\nu)$,  similar to \eqref{eq:estimate}, it suffices to prove that
\beq \label{eq:requires}
\sum_{n \ge N} |(\ip{f}{\gamma_\mu})_{(n)}|^2 \leq C \cdot N^{-2} \cdot D(N) \qquad \text{as } N \rightarrow \infty
\eeq
with $N^{-\epsilon} \cdot D(N) \to 0$ as $N \to 0$ for all $\epsilon > 0$. We remark that in the following
the constants might change, by abuse of notation, we however always coin them $C$.

First recall that, by \cite[Thm 1.1]{GL07a}, the shearlet frame $(\sigma_\eta)_\eta$ achieves the rate
\beq \label{eq:req1}
\sup_{g \in \cE^2(\nu)} |(\ip{g}{\sigma_\eta})_{(n)}| \leq C \cdot n^{-3/2} \cdot (\log n)^{3/2} \qquad \text{for each } n.
\eeq
Now let $\epsilon > 0$, and choose $p = 2/3 + \epsilon$. Then, by \eqref{eq:req1},
\[
\sup_{g \in \cE^2(\nu)} \norm{(\ip{g}{\sigma_\eta}_\eta)}_p^p \le \sup_{g \in \cE^2(\nu)} C \cdot \sum_n (n^{-3/2} \cdot (\log n)^{3/2})^{2/3 + \epsilon}
\le C.
\]
By Theorem \ref{theo:main2}, this implies that $\sup_{g \in \cE^2(\nu)} \norm{(\ip{g}{\gamma_\mu}_\mu)}_p \le C$.
Hence, for each $n$,
\[
\sup_{g \in \cE^2(\nu)} |(\ip{g}{\gamma_\mu})_{(n)}| \le C \cdot n^{-1/p},
\]
and therefore
\[
\sum_{n \ge N} |(\ip{f}{\gamma_\mu})_{(n)}|^2 \leq C \cdot N^{-2/p+1} = C \cdot N^{-2} \cdot N^{\frac{9\epsilon}{2+3\epsilon}}.
\]
In the definition of $p$ the variable $\epsilon$ can be chosen arbitrarily small, which implies \eqref{eq:requires},
and the theorem is proved. \qed

%
%
%
%
%



\subsection{Geometric Separation}

Natural images are typically composed of morphologically distinct features; an example being
spines (pointlike structures) and dendrites (curvelike structures) in neurobiological
imaging. One goal is to automatically extract those
components for separate analysis. In \cite{DK08}, the author, joint with Donoho, studied the
situation of images composed of point- and curvelike structures, for which they introduced
models by
\beq \label{eq:CandP}
   \cP = \sum_{i=1}^P |x - x_i|^{-3/2}
   \quad \mbox{and} \quad
   \cC = \int \delta_{\tau(t)}dt, \quad \mbox{ with } \tau: [0,1] \mapsto \bR^2\mbox{ a closed curve,}
\eeq
respectively. The {\em Geometric Separation Problem} now consists in extracting $\cP$ and $\cC$ from knowledge
of $f$ given by
\[
f = \cP + \cC.
\]
In \cite{DK08}, a particular decomposition technique based on $\ell_1$ minimization was employed which
required suitably chosen overcomplete systems which sparsify the different components. Using the tight frame
of radial wavelets for the pointlike structures and the tight frame of curvelets for the curvelike structures,
asymptotically arbitrarily precise separation was proven.

Using the results on sparsity equivalence derived in this paper, we can now prove that a different pair
of representation systems can be utilized for this Geometric Separation Problem, which is more suitable
for a digital realization: orthonormal separable Meyer wavelets and shearlets. In contrast to the pair
considered before, surprisingly, now one system even forms an orthonormal basis.
%

\begin{figure}[ht]
\begin{center}
\includegraphics[height=1.75in]{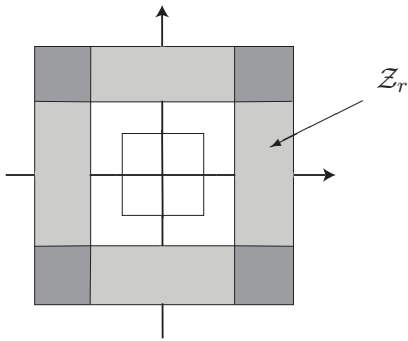}
\put(10,90){\vector(-2,-1){35}}
\put(15,95){$\cZ_r$}
\end{center}
\caption{The tiling of the frequency domain induced by orthonormal separable Meyer wavelets.}
\label{fig:ONWavelets}
\end{figure}

For the reader's convenience, we first briefly recall the definition of orthonormal separable Meyer wavelets.
Let $W \in L^2(\bR)$ denote the Fourier transform of the Meyer wavelet and $\phi \in
L^2(\bR)$ the associated scaling function. Letting $W^h \in L^2(\bR^2)$, $h=1,2,3$ be defined by
\[
W^1(\xi)=\hat{\phi}(\xi_1) W(\xi_2),\quad W^2(\xi)=W(\xi_1) \hat{\phi}(\xi_2)\quad \mbox{and} \quad
W^3(\xi)=W(\xi_1) W(\xi_2),
\]
the orthonormal separable Meyer wavelets at scale $j$ and spatial position $n$ are defined by their Fourier transforms
\[
\hat{\psi}_{\nu}(\xi) = 2^{-j} W^h(\xi/2^j) e^{i\ip{n}{\xi}/2^j},
\]
where $\nu=(h,j,n)$ index type of mother function, scale, and position. This system forms an orthonormal basis
for $L^2(\bR^2)$. For each $j$, the functions $\hat{\psi}_{\nu}$ are supported on the corona $\cZ_{2^{j+1}\pi/3}$, where
\[
\cZ_r = \{\xi \in \bR^2 : r \le \norm{\xi}_\infty \le 4 \cdot r\}
\]
(see Figure \ref{fig:ONWavelets}). For more details we refer to \cite{Mal98}.

Shearlets $\{\sigma_\eta\}_\eta$, where $\eta = (j,k,m,\iota)$ indexes scale,
orientation, position, and cone, were already defined in Subsection \ref{subsec:shearlets}, but
to match them with Meyer wavelets, we now choose $W$ to be the Fourier transform of the Meyer wavelet.
We wish to draw the reader's attention to the fact that the supports
of orthonormal separable Meyer wavelets match perfectly with the
supports of shearlets. In fact, for each scale $j$, the Fourier transforms of
the elements of both systems are supported on $\cZ_{2^{j+1}\pi/3}$.


We next construct a family of filters $F_j$ with transfer functions
\[
    \hat{F}_j(\xi) = W(\|\xi\|_\infty/2^{j}), \qquad \xi \in \bR^2
\]
leading to a decomposition of a function $g$ into
functions $g_j = F_j \star f$ defined on the frequency corona
$Z_{2^{j+1}\pi/3}$ equipped with the reconstruction formula $g =
\sum_{j} F_j \star g_j$. Let $\cF_j$ denote  the range of the
operator of convolution with $F_j$. Then shearlets at level $j'$ are
orthogonal to $\cF_j$ unless $|j' -j | \leq 1$.  Similarly,
orthonormal separable Meyer wavelets at level $j'$ are orthogonal to
$\cF_j$ unless $|j' -j | \leq 1$. The proofs of these two claims use
precisely the same arguments as the corresponding result in
\cite{DK08}, wherefore we omit them.

We can now formulate the corresponding Component Separation Problem (CSep). For the sake of brevity, we let
$\Theta_j$ denote the indices $\nu = (h,j,n)$ of orthonormal separable Meyer wavelets at level $j$, and let
$\Theta_j^\pm = \Theta_{j-1} \cup \Theta_j \cup \Theta_{j+1}$. Likewise, we let $\Sigma_j$ denote the indices
$\eta = (j,k,m,\iota)$ of shearlets at level $j$, and let $\Sigma_j^\pm = \Sigma_{j-1} \cup \Sigma_j \cup \Sigma_{j+1}$.
Further, we denote the filtered composed image $f$ and the filtered point and curvilinear part $\cP$ and $\cC$
(cf. \eqref{eq:CandP}) by
\[
f_j = F_j \star f = F_j \star (\cP + \cC) = \cP_j + \cC_j.
\]
Then we can formulate the Component Separation Problem as the following $\ell_1$ minimization problem:
\[
(\mbox{\sc CSep})  \qquad   (W_j,S_j) = \mbox{ argmin } \|(\langle W_j , \psi_\nu \rangle)_\nu\|_1 + \|(\langle S_j , \sigma_\eta \rangle)_\eta\|_1
\quad \mbox{subject to }  f_j = W_j + S_j.
\]
We claim that the considered pair of representation systems leads to
asymptotically perfect separation in the sense of the following theorem. Before stating the result, we wish
to remark that the proof draws from various definitions and lemmata from \cite{DK08}, wherefore we decided that
for the sake of brevity -- this being mostly an application of our main result in this paper -- we only present
the road map of its proof.

\begin{theorem}
\label{maintheorem2}
Let $(W_j,S_j)$ denote the solution of $(\mbox{\sc CSep})$. Then, we have
\[
\frac{ \| W_j - \cP_j \|_2 + \| S_j - \cC_j \|_2 }{\| \cP_j\|_2 +
\|\cC_j\|_2 } \goto 0, \qquad j \goto \infty
\]
\end{theorem}

\noindent
{\bf Proof.}
The proof presented in \cite{DK08} uses as one main idea the following estimate for each scale $j$: Let $\cS_{1,j}$ and
$\cS_{2,j}$ be sets of `significant coefficients' of wavelets and curvelets, respectively,
let $\delta_j$ be the sparse approximation error given by
\[
\sum_{\nu \in \cS_{1,j}^c} |\langle W_j , \psi_\nu \rangle| +  \sum_{\eta \in \cS_{2,j}^c} |\langle S_j , \sigma_\eta \rangle| \le \delta_j,
\]
and let $(\mu_c)_j$ be the cluster coherence defined as
\[
(\mu_c)_j = \max\Big\{\max_\eta \sum_{\nu \in \cS_{1,j}} \absip{\psi_\nu}{\sigma_\eta}, \max_\nu \sum_{\eta \in \cS_{2,j}} \absip{\psi_\nu}{\sigma_\eta}\Big\}.
\]
Then \cite[Prop. 2.1]{DK08} applied to each filtered $f_j$ implies
\[
\| W_j - \cP_j \|_2 + \| S_j - \cC_j \|_2 \le \frac{2\delta_j}{1-2(\mu_c)_j}.
\]

Thus, the key step in \cite{DK08} was the construction of clusters $\cS_{1,j}$ and $\cS_{2,j}$  having {\em both} of the following
two properties: (i) asymptotically negligible cluster coherences:
\[
   (\mu_c)_j  \goto 0, \qquad j \goto \infty,
\]
and (ii) asymptotically negligible cluster approximation errors:
\[
  \delta_j = o(\|\cP_j\|_2 + \|\cC_j\|_2) , \qquad j \goto \infty.
\]

The same steps with very similar argumentations can be performed for the pair wavelets-shearlets if adapted
clusters $\cS_{1,j}$ and $\cS_{2,j}$ are defined by applying the following two key observations:
\begin{itemize}
\item It was shown in Theorem \ref{theo:main2} that shearlets and curvelets are
sparsity equivalent; more precisely, there exists a sparse matrix
$N_1$, say, which satisfies
\[
(\ip{\sigma_\eta}{g})_\eta = N_1
(\ip{\gamma_\mu}{g})_{\mu}
\]
for any distribution $g$.
\item Orthonormal separable Meyer wavelets and radial wavelets are likewise
sparsity equivalent, i.e., there exists a sparse matrix $N_2$, say,
which satisfies
\[
(\ip{\psi_\nu}{g})_\nu = N_2 (\ip{\psi_\lambda}{g})_{\lambda}
\]
for any distribution $g$.
\end{itemize}

A second ingredient are estimates for inner products between wavelets and shearlets within the frames, but also across.
For this, the paralleling lemma to \cite[Lem. 3.3]{DK08} -- with a very similar proof -- is essential:
\begin{lemma}
For each $N = 1,2, \dots$ there is a constant $c_N$ so that
\[
   \absip{\psi_\nu}{\psi_\lambda} \leq c_N \cdot  1_{\{|j-j'| < 2\}}  \cdot
    \langle |n - n'| \rangle^{-N},
   \qquad \forall \nu = (h,j,n) \; \forall \lambda = (h',j',n').
\]
\end{lemma}

As already remarked before, we will not lay out the precise details of the complete proof, since the arguments in the very lengthy and
technical proof from \cite{DK08} just need to be adapted in a straightforward manner to the sets of
significant coefficients now based on the choice for orthonormal wavelets and shearlets. We then derive Theorem
\ref{maintheorem2}, thus perfect separation using orthonormal separable Meyer wavelets and shearlets.
\qed


\section{Extensions and General Viewpoint}
\label{sec:extensions}

So far we focused entirely on a very special situation showing sparsity equivalence between
curvelets and shearlets. Our goal was to show that for this exemplary situation sparsity equivalence
can be established, provides insight into the relation between these systems, and lead automatically
to novel results on sparse expansions of those two anisotropic systems.

This is however just the `tip of the iceberg': the main results in this paper are susceptible of very extensive
generalizations and extensions.

\bitem
\item {\em Curvelets and Shearlets}. It is conceivable that a similar statement as Theorem \ref{theo:main2} is
provable for first generation curvelets as also for the new class of compactly supported shearlets. It should
though be mentioned that the compactly supported shearlet frames introduced so far are not tight frames,
hence the framework developed in this paper needs to be extended to pairs of general frames.
\item {\em Other Systems}. The analysis of sparsity equivalence of curvelets and shearlets we drove here can
and should be applied to other pairs of systems. Ideally, novelly introduced systems could be compared to
a system whose sparse approximation properties are already very well understood.
\item {\em Systems with Continuous Parameters}. Certainly, we can also ask about similar sparsity
properties for systems with continuous parameters. This however requires a different sparsity model, where one
conceivable path would be to compare resolution of wavefront set behavior in the sense of \cite{CD05a,KL07}.
\item {\em Weighted Norms}. When aiming at transferring results such as sparse decompositions of curvilinear
integrals \cite{CD00a} or sparse decompositions of the Radon transform \cite{CD00b}, the framework needs to be
generalized to weighted $\ell_p$ norms. Also the analysis of associated approximation spaces requires this
extension, since, for instance, the norm associated with the curvelet spaces introduced in \cite[p.~67]{BN07}
is precisely a weighted mixed $\ell_{p,q}$ norm of the coefficient sequence.
\eitem

\end{document}